\documentclass[10pt]{article}
\usepackage{isolatin1}
\usepackage{multirow}
\usepackage{calc}
\usepackage{epsfig,amsfonts,amsthm,latexsym,amsmath}
\usepackage{hyperref}
\usepackage{color}
\usepackage{graphicx}
\usepackage{caption}
\usepackage{subcaption}
\usepackage{tikz}

\title{\normalsize\bf%
SOLVING IRREGULAR STRIP PACKING PROBLEMS WITH FREE ROTATIONS USING SEPARATION LINES
}

\author{%
Jeinny Peralta$^{1,2,*}$,\ \ Marina Andretta$^{1}$ \ and \ Jos\'e Fernando Oliveira$^{3}$
}

\begin{document}

\date{}

\maketitle

\vspace{-20pt}
\begin{center}
{\footnotesize 
*Corresponding author\\
$^1$Universidade de S\~ao Paulo, Av.Trabalhador S\~ao-carlense, 400, 13566-590, S\~ao Carlos - SP, Brasil \\
$^2$Universidad de la Costa, Calle 58 \# 55 - 66, Barranquilla, Colombia\\
$^3$Universidade do Porto, Rua Dr. Roberto Frias, 4200-590, Porto, Portugal\\
E-mails: jperalta10@cuc.edu.co / andretta@icmc.usp.br /  jfo@fe.up.pt  
}\end{center}

\bigskip
\noindent
{\small{\bf ABSTRACT.}
Solving nesting problems or irregular strip packing problems is to position polygons in a fixed width and unlimited length strip, obeying polygon integrity containment constraints and non-overlapping constraints, in order to minimize the used length of the strip.  To ensure non-overlapping, we used separation lines.  A straight line is a separation line if given two polygons, all vertices of one of the polygons are on one side of the line or on the line, and all vertices of the other polygon are on the other side of the line or on the line. Since we are considering free rotations of the polygons and separation lines, the mathematical model of the studied problem is nonlinear. Therefore, we use the nonlinear programming solver IPOPT (an algorithm of interior points type), which is part of COIN-OR. Computational tests were run using established benchmark instances and the results were compared with the ones obtained with other methodologies in the literature that use free rotation.  

}


\baselineskip=\normalbaselineskip

\section{Introduction}\label{sec:1}

Irregular strip packing problems have a great relevance in the production processes, such as garment manufacturing and furniture making. In irregular strip packing, smaller irregular pieces (in our case polygons) must be positioned into a big piece (in our case the strip), minimizing the used length of the strip.  The main constraint in irregular strip packing problems is the non-overlapping between pieces, but it is very complex for a computational program to determine if two pieces are overlapping, touching or separated.  In the literature there are tools for solving this issue (\cite{bennell2008}), among these the raster methods, direct trigonometry, no-fit polygon, and phi-function.  In raster methods the strip is always divided into discrete areas and coding schemes are used.  In the coding schemes used by Oliveira and Ferreira, and Segenreich and Braga in \cite{b1} and \cite{SB}, respectively, the empty cells belonging to the division of the strip are encoded by zero and numbers equal to or greater than one are used to encode a piece; so, to check the non-overlap in raster methods is only a matter of checking the grid cells. There are several tools that use trigonometry to deal with the non-overlapping. In \cite{jones2013} circles inscribed are used to relax the non-overlapping constraints, replacing them with non-overlapping constraints of circles inscribed.  In the remaining tools that use trigonometry, the pieces are represented by polygons, thereby, the non-intersection of the edges of polygons must be checked to check the non-overlap.
In phi-function, the pieces are represented by the union or intersection of primary objects, that is, circles, rectangles, regular polygons, convex polygons, and the complement of these forms.  This tool was designed and implemented in \cite{stoyan2001, stoyan2004, bennell2010, stoyan15}.  The phi-function is a mathematical expression that represents the relative position of two pieces.  Specifically, the phi-function value is greater than zero if the pieces are separated; equal to zero if their borders are touching; and smaller than zero if overlapping each other.  In this article, we use two tools to ensure non-overlapping of pieces, direct trigonometry and no-fit polygon.  In the model, we represent the pieces by polygons and use separation lines, that is, we use trigonometry.  For the resolution of the modeled problem we need a starting point; for the construction of this point, we use no-fit polygon.  A straight line is a separation line if all vertices of a polygon are on one side of the line or on the line, and all vertices of the other polygon are on the other side of the line or on the line.  In no-fit polygon the pieces are also represented by polygons.  The no-fit polygon is a polygon resulting from the two polygons that are being compared.  One of the advantages of this method is that the generation of these no-fit polygons is done only once, in a pre-processing phase, but a big disadvantage is that the no-fit polygons are dependent on the orientation of the polygons, and have to be generated for all possible orientations.  

Several solving techniques to these problems that deal with irregular pieces have been developed, based predominantly on heuristics and metaheuristics (\cite{b1,b6,b7,b3,b8,b9}). The heuristics used for solving these problems can work with partial solutions, constructing the final layout piece by piece (constructive heuristics), or complete solutions, in which changes are done in order to find improvements.  Exact algorithms of mixed integer programming that ensure finding the optimal solution were also developed (\cite{fischetti,alvarez-valdes,toledo}); however, in these algorithms, the runtime increases dramatically with the increase of the quantity of objects used in the problem. In these techniques and algorithms, free rotations are not allowed.

Additionally in the literature, we also find a visual system for packing problem of irregular pieces with free rotation into a rectangular board that aims to minimize the waste, \cite{Liao16}.  This algorithm is based on a method of physics, the rubber band physics movement.  

Nonlinear programming models have also been proposed for representing the irregular packing problem, such as \cite{chernov2010, jones2013, rochaetal,kallrath2009,stoyan15,Liao16}. In all these models free rotation of the pieces is allowed. 
In \cite{chernov2010} a model for a strip packing problem was presented.  In this model phi-functions are used to ensure non-overlapping of the pieces.  The pieces in this paper are phi-objects, which are 2D and 3D objects of very general type.  To solve the problem Chernov, Stoyan and Romanova applied a modification of the Zoutendijk method of feasible directions \cite{Z1960,Z1970} combined with the concept of $\epsilon$-active inequalities \cite{Stoyan2008}.  
In \cite{jones2013} the pieces and the shapes can be arbitrary no convex polygons and to solve the problem Jones used three solvers: Branch\&Reduce Optimization Navigator (BARON) \cite{ts05,sahinidis:baron:14.3.1}; LindoGlobal from Lindo Systems, Inc., which is part of the GAMS 22.5 distributions; and GloMIQO \cite{MIQO}.  
In \cite{rochaetal}, like in \cite{chernov2010}, a model for an irregular strip packing problem was presented.  In this, the resolution of the problem is divided in two phases: big pieces are compacted in a first phase, while in a second phase, the remaining small pieces are placed between the big pieces. In their experiments, Rocha et al. used instances where the pieces are convex and no convex irregular polygons.  The representation of these polygons was done by circle covering, and they used the nonlinear solver ALGENCAN \cite{andreani2007,andreani2008}.
In \cite{kallrath2009} a model for two cases of packing pieces was developed. In the first case, the objective is to pack the pieces in such a way as to minimize the area of the design rectangle. In the second case, the objective is to pack the pieces on stocked rectangles of known geometric dimensions.  Separation lines are used to ensure non-overlapping.  In their work the pieces are circles, rectangles, and convex polygons; and to solve the problem Kallrath used BARON \cite{ts05,sahinidis:baron:14.3.1} and LindoGlobal; he performed an experiment with only two polygons and found a feasible solution for it, in which LindoGlobal proved global optimality in 40 min, but BARON did not increase the lower bound at all.  With more than two polygons this technique has difficulty finding optimal solution to the problem. 
In \cite{stoyan15} was provided a nonlinear programming model that employs ready-to-use phi-functions.   In their paper, the pieces are bounded by circular arcs and/or line segments, and two types of container are considered, rectangular and circular. To solve the problem, Stoyan, Pankratov and Romanova developed a compaction algorithm to search for local optimal solutions, which is performed by IPOPT (an algorithm of interior points type, \cite{ipopt}), which is part of the COIN-OR. 

In this paper we propose a nonlinear mathematical model for an irregular strip packing problem which deals only with polygons which may be convex or no convex, and that can rotate freely.  In the model, to ensure non-overlapping, we use direct trigonometry, in particular separation lines, a similar technique to that used in \cite{kallrath2009}, but with a significantly lower number of variables, allowing us to obtain good solutions for larger instances in reasonable execution times.  As said before, a straight line is a separation line if given two polygons, all vertices of one of the polygons are on one side of the line or on the line, and all vertices of the other polygon are on the other side of the line or on the line. Like the polygons, the separation lines also can rotate freely.  We use a code for nonlinear programming to solve the problem, IPOPT \cite{ipopt}, which depends substantially on a starting point.  We present a way of calculating starting points.  

This paper is organized as followed.  In the next section a model of an irregular strip packing problem that considers free rotations is presented.  The modeled problem, the polygons representation used in the model, and the tool used to ensure non-overlapping, are described also in this section.  In Section \ref{resultados}, the parameters of the algorithm used for solving the problem are presented, as well as the numerical results obtained when performing tests with known benchmark instances.  At the end of our paper we present some conclusions in Section \ref{conclusion}.

\section{A model for an irregular strip packing problem}\label{sec:2}

The irregular strip packing problem studied in this paper consists of placing $n$ irregular polygons, which can rotate freely, in a fixed width and unlimited length strip, obeying polygon integrity containment constraints and non-overlapping constraints, in order to minimize the used length of the strip.  We propose a nonlinear mathematical model for this irregular strip packing problem. 

We now explain how the polygons are represented in the model (Section \ref{representation}), as well as the tool used to ensure non-overlapping of the polygons (Section \ref{separator lines}), and then, introduce the model (Section \ref{modelo}).

\subsection{Representation of polygons in the model}
\label{representation}

Here, we describe the representation of the polygons.  Remember that the polygons may be convex or no convex. If a polygon is convex, it is represented by their vertices, as follows:

$$P_i=[(x_{i}^1,y_{i}^1),(x_{i}^2,y_{i}^2),\ldots,(x_{i}^{v_i},y_{i}^{v_i})],$$
being $v_i$ the number of vertices of the polygon $P_i$.  

\begin{figure}[!h]
\centering
\begin{tikzpicture}[thick,x=5,y=5]
\draw[draw=, fill=, fill opacity=0.35](50.000000,20.689655)--(55.172414,15.517241)--(58.620690,15.517241)--(53.448276,20.689655)--cycle;
\node[scale=0.8]  at (54.800000,17.689655){$P_{1_1}$};
\draw[draw=, fill=, fill opacity=0.35](46.551724,25.862069)--(50.000000,20.689655)--(53.448276,20.689655)--(50.000000,25.862069)--cycle;
\node[scale=0.8]  at (50.400000,22.889655){$P_{1_2}$};
\draw[draw=, fill=, fill opacity=0.35](46.551724,31.034483)--(46.551724,25.862069)--(50.000000,25.862069)--(50.000000,31.034483)--cycle;
\node[scale=0.8]  at (48.400000,28.689655){$P_{1_3}$};
\draw[draw=, fill=, fill opacity=0.35](50.000000,36.206897)--(46.551724,31.034483)--(50.000000,31.034483)--(53.448276,36.206897)--cycle;
\node[scale=0.8]  at (49.700000,33.389655){$P_{1_4}$};
\draw[draw=, fill=, fill opacity=0.35](55.172414,41.379310)--(50.000000,36.206897)--(53.448276,36.206897)--(58.620690,41.379310)--cycle;
\node[scale=0.8]  at (53.700000,38.289655){$P_{1_5}$};
\end{tikzpicture}
\caption{Partition of a no convex polygon.}
\label{partition}
\end{figure}
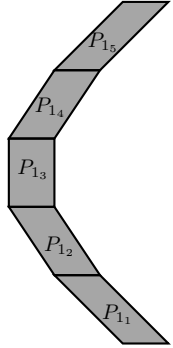

 If a polygon is no convex, it can be decomposed in convex polygons, as follows: 
$$P_i=[P_{i_1}, P_{i_2},\ldots, P_{i_{p_i}}],$$ 
being $p_i$ the number of convex polygons belonging to the partition of the no convex polygon $P_i$, see Figure~\ref{partition}. 


The coordinates of a vertex belonging to a partition of a no convex polygon $P_i$ are given by $$(x_{i_j}^l,y_{i_j}^l),$$ with $j=1,\ldots,p_i$ and $l=1,\ldots,v_{i_j}$, being $v_{i_j}$ the number of vertices of the convex polygon $P_{i_{j}}$.

We can deal with the problem with $n$ no convex polygons in the same way that we deal the problem with $N$ convex polygons, with $N=\sum_{i=1}^{n}p_i$, we just have to ensure that the translations and rotations are the same for all polygons belonging to the partition of a no convex polygon.

The reference point is used for representing a polygon which has undergone translations and/or rotations, since we can write all other vertices of the polygon in terms of this point, as can be seen in Section \ref{vertices}.    

\begin{figure}[!h]
\centering
\begin{tikzpicture}[thick,x=1.5,y=1.5]
\draw[draw=, fill=, fill opacity=0.35](119.714286,43.224490)--(100.000000,44.367347)--(100.000000,1.755102)--(119.714286,2.897959)--cycle;
\draw[draw=, fill=, fill opacity=0.35](144.591837,41.265306)--(140.469388,46.122449)--(119.714286,43.224490)--(119.714286,2.897959)--(140.469388,0.000000)--(144.591837,4.857143)--cycle;
\draw[draw=, fill=, fill opacity=0.35](144.591837,4.857143)--(155.795918,4.428571)--(161.224490,15.653061)--(157.530612,18.367347)--cycle;
\draw[draw=, fill=, fill opacity=0.35](144.591837,4.857143)--(157.530612,18.367347)--(157.530612,27.755102)--(144.591837,41.265306)--cycle;
\draw[draw=, fill=, fill opacity=0.35](157.530612,27.755102)--(161.224490,30.469388)--(155.795918,41.693878)--(144.591837,41.265306)--cycle;
\draw[ -> ](100,1.755102)--(100,60.000000)node[scale=0.7] [left] {$y$};
\draw[->](100,1.755102)--(170,1.755102)node[scale=0.7] [below]{$x$};
\fill (100,1.755102) circle(1.5pt);
\draw (119.714286,2.897959) circle(1.5pt);
\draw (144.591837,4.857143) circle(1.5pt);
\draw (157.530612,27.755102) circle(1.5pt);
\node[scale=0.7] at (100.000000,-0.55555) {Reference Point};
\end{tikzpicture}
\caption{Reference point of a no convex polygon}
\label{reference}
\end{figure}
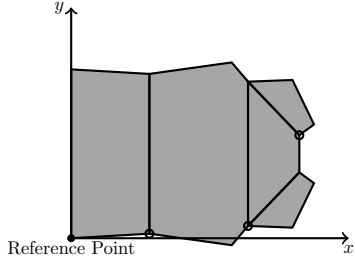

We choose a vertex nearer to the origin as the reference point of a polygon, see Figure~\ref{reference}. 


\subsubsection{Vertices of polygons in general form}
\label{vertices}

Henceforth, we will use the following notation: $(x_{i}^l,y_{i}^l)$ are the coordinates of a vertex of a polygon $P_i$ in the original position, and $(\bar{x}_{i}^l,\bar{y}_{i}^l)$ are the coordinates of a vertex of a polygon $P_i$ which has undergone translations and/or rotations.


In the representation of the polygons used in the model, all the vertices are translated, so that the reference point (the first vertex) is located at the origin, $(x_i^1,y_i^1)=(0,0)$.   Note that the reference point will be transferred after rotation and translation, to the point $(t_{i_x},t_{i_y})$, that is, $(\bar{x}_{i}^1,\bar{y}_{i}^1)=(t_{i_x},t_{i_y})$, with $t_{i_x}$ and $t_{i_y}$ the values of the translation parameters of the polygon $P_i$.  Thus, if we translate and rotate a polygon $P_i$ around his reference point, the coordinates of the vertices in general form are given by:    

\begin{equation}
 (\bar{x}_{i}^l,\bar{y}_{i}^l)=(x_{i}^l\cos\theta_i - y_{i}^l\sin\theta_i + \bar{x}_{i}^1, x_{i}^l\sin\theta_i + y_{i}^l\cos\theta_i + \bar{y}_{i}^1),
\end{equation} 

\noindent being $\theta_i$ the angle of rotation of the polygon $P_i$. We consider that positive angles represent rotation in the counterclockwise direction.


When we deal with no convex polygons, we make sure that the translations and rotations are the same for all polygons belonging to the partition of the no convex polygon.

\subsection{Separation lines}
\label{separator lines}

We use separation lines to ensure non-overlapping. A straight line given by the equation $y = c_{i,r}x + d_{i,r}$ separates two polygons $P_i$ and $P_r$ if either
\begin{equation}\label{r1} 
\left\{
\begin{aligned}
& y_{i}^l - c_{i,r}x_{i}^l - d_{i,r} \leq 0, \quad \forall l\in\{1,\ldots,v_i\},\\ 
& y_{r}^l - c_{i,r}x_{r}^l - d_{i,r} \geq 0, \quad \forall l\in\{1,\ldots,v_r\},
\end{aligned}
\right.
\end{equation} 
or
\begin{equation}\label{r2}
\left\{
\begin{aligned}
& y_{i}^l - c_{i,r}x_{i}^l - d_{i,r} \geq 0, \quad \forall l\in\{1,\ldots,v_i\},\\ 
& y_{r}^l - c_{i,r}x_{r}^l - d_{i,r} \leq 0, \quad \forall l\in\{1,\ldots,v_r\}.
\end{aligned}
\right.
\end{equation}\\  

\begin{figure}[!hb]
\centering
\begin{tikzpicture}[thick,x=1.5,y=1.5]
\draw[draw=, fill=, fill opacity=0.35](119.714286,43.224490)--(100.000000,44.367347)--(100.000000,1.755102)--(119.714286,2.897959)--cycle;
\draw[draw=, fill=, fill opacity=0.35](164.591837,41.265306)--(160.469388,46.122449)--(129.714286,43.224490)--(129.714286,2.897959)--(160.469388,0.000000)--(164.591837,4.857143)--cycle;
\draw[dashed][line width=0.35mm, red](119.714286,50.224490)--(119.714286,-5.897959);
\draw[ -> ](90,-10)--(90,60.000000)node[scale=0.75] [left] {$y$};;
\draw[->](90,-10)--(170,-10)node[scale=0.75] [below]{$x$};
\end{tikzpicture}
\caption{Separation line of two convex polygons.}
\label{figretas}
\end{figure}
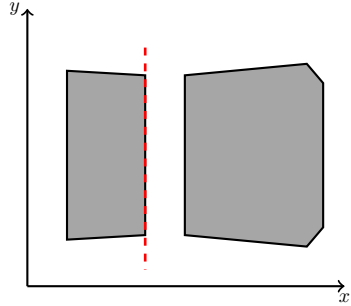

That is, a straight line is a separation line if all vertices of a polygon are on one side of the line or on the line, and all vertices of the other polygon are on the other side of the line or on the line, see Figure~\ref{figretas}. 


When we are dealing with no convex polygons, we must have lines separating each pair of polygons $P_{i_j}$, $P_{r_s}$, belonging to the partition of polygons $P_i$, $P_r$, respectively, with $i \neq r$ (that is, we do not have lines separating the polygons belonging to the partition of a no convex polygon), $j\in {1,\ldots,p_i}$, and $s\in {1,\ldots,p_r}$./


Next, we present the general form of a separation line passing through a side of one of the two polygons, that is, which passes through two vertices of a polygon ($P_{i_j}$ or $P_{r_s}$), let's say $(x_k , y_k)$ and $(x_{k+1}, y_{k+1})$.

\begin{equation}\label{r,s}
\begin{aligned}
(y - y_k)(x_{k+1} - x_k) - (y_{k+1} - y_k)(x - x_k) = 0.
\end{aligned}
\end{equation}

Like polygons, the separation lines can also rotate and translate, as long as they remain being separation lines. When rotating and translating a separation line, we rewrite the point $(\bar{x}_{k+1}, \bar{y}_{k+1})$, in function of $(\bar{x}_k,\bar{y}_k)$, which we will call from now on, reference point of the separation line; therefore

\begin{equation*}
\begin{aligned}
 (\bar{x}_{k+1}, \bar{y}_{k+1})=((x_{k+1} - x_{k})\cos\alpha_{i_j,r_s} + &(y_{k} - y_{k+1})\sin\alpha_{i_j,r_s} + \bar{x}_k, \\ 
(x_{k+1} - x_{k})\sin\alpha_{i_j,r_s} + &(y_{k+1} - y_k)\cos\alpha_{i_j,r_s} + \bar{y}_k),
\end{aligned}
\end{equation*}

\noindent being $\alpha_{i_j,r_s}$ the angle of rotation of the separation line of the polygons $P_{i_j}$ and $P_{r_s}$.

Next, we rewrite the separation line equation:

\begin{equation}\label{r,s1}
\begin{aligned}
(y - \bar{y}_k)&[(x_{k+1} - x_{k})\cos\alpha_{i_j,r_s} + (y_{k} - y_{k+1})\sin\alpha_{i_j,r_s}] - \\ 
(x - \bar{x}_k)&[(x_{k+1} - x_{k})\sin\alpha_{i_j,r_s} + (y_{k+1} - y_k)\cos\alpha_{i_j,r_s}] = 0.
\end{aligned}
\end{equation}\\  

\subsection{Model for an irregular strip packing problem with $n$ polygons considering free rotations}\label{modelo}

Because we want to position $n$ polygons in a fixed width and unlimited length strip in order to minimize the used length of the strip, the objective function is given by:
 
\begin{equation}\label{z}
z=\max\{\bar{x}_{i_{l}}\}, \quad l=1,2,\ldots,v_i \text{ and } i=1,2,\ldots,n,
\end{equation}

\noindent being $v_i$ the number of vertices of the polygon $P_i$.

Without loss of generality, to ensure non-overlapping, we used in our model the set of constraints \ref{r1}, for each pair of polygons  $P_{i_j}$ and $P_{r_s}$, with  $$c_{i_j,r_s}=\frac{(x_{k+1} - x_{k})\sin\alpha_{i_j,r_s} + (y_{k+1} - y_k)\cos\alpha_{i_j,r_s}}{(x_{k+1} - x_{k})\cos\alpha_{i_j,r_s} + (y_{k} - y_{k+1})\sin\alpha_{i_j,r_s}}$$ and $$d_{i_j,r_s}=\bar{y}_k - c_{i_j,r_s}\bar{x}_k,$$

\noindent in which $(x_k , y_k)$ and $(x_{k+1}, y_{k+1})$ are the two vertices of one of the polygons whereby passes the separation line and $\alpha_{i_j,r_s}$ is the rotation angle of the straight line that separates $P_{i_j}$ from $P_{r_s}$.  Note that the reference point of the separation line, $(\bar{x}_k,\bar{y}_k)$ are the values of the translations parameters, that from now on, we will write $(\bar{x}_{i_j,r_s}, \bar{y}_{i_j,r_s})$, for a straight line that is separating  $P_{i_j}$ from $P_{r_s}$.

Let $(z,q_1,q_2,...,q_n,\bar{q}_1,\bar{q}_2,...,\bar{q}_Q)$ be the vector of all variables in our model, being $z$ the length of the strip defined in \ref{z}, $q_i$ the variables referring to the polygon $P_i$, $q_i=(\bar{x}_i^1,\bar{y}_i^1,\theta_i)$, $i=1,..,n$ and  $\bar{q}_\ell$ the variables referring to the line that separates polygon $P_{i_j}$ from polygon $P_{r_s}$, $\bar{q}_\ell=(\bar{x}_{i_j,r_s}, \bar{y}_{i_j,r_s},\alpha_{i_j,r_s})$, $\ell=1,...,Q$ and $Q=\sum_{i=1}^{n-1}p_i(N-\sum_{k=1}^ip_i)$ the number of separation lines. Let $e$ be the width of the strip in which the polygons are to be placed.  A general model for our problem is given by:\\

\begin{subequations}
\begin{alignat}{4}\nonumber
\text{Minimize } & z\\\label{ec1}
\text{subject to } & 0\leq \bar{y}_{i_j}^l\leq e, &i&=1,\ldots,n,\\\nonumber
                  &                       &j&=1,\ldots,p_i,\\\nonumber
                  &                       &l&=1,\ldots,v_{i_j},\\\label{ec2}
                  & 0\leq \bar{x}_{i_j}^l\leq z, &i&=1,\ldots,n,\\\nonumber
                  &                       &j&=1,\ldots,p_i,\\\nonumber
                  &                       &l&=1,\ldots,v_{i_j}, \\\label{ec3}
                  & \bar{y}_{i_j}^l - c_{i_j,r_s}\bar{x}_{i_j}^l - d_{i_j,r_s} \leq 0, &  i&=1,\ldots,n,\\ \nonumber
                  &                                          &  r&=1,\ldots,n,\\\nonumber
\nonumber                  &                                          &  i&\neq r,\\
\nonumber                  &                                          &j&=1,\ldots,p_i,\\
\nonumber                  &                       &s&=1,\ldots,p_r,\\  
\nonumber                  &                        &l&=1,\ldots,v_{i_j},\\ \label{ec4}  
                  &\bar{y}_{r_s}^l - c_{i_j,r_s}\bar{x}_{r_s}^l - d_{i_j,r_s} \geq 0,&  i&=1,\ldots,n,\\   
\nonumber                  &                                          &  r&=1,\ldots,n,\\
\nonumber                  &                                          &  i&\neq r,\\
\nonumber                  &                      &j&=1,\ldots,p_i,\\                    
\nonumber                 &                       &s&=1,\ldots,p_r.\\    
\nonumber                  &                      &l&=1,\ldots,v_{r_s}, \\ \nonumber
\end{alignat}
\end{subequations} \\

Remembering that the vertices of a translated and rotated polygon $P_{i_j}$ are given by $(\bar{x}_{i_j}^l,\bar{y}_{i_j}^l)$, for  $i=1,\ldots,n$, $j=1,\ldots,p_i$ and $l=1,\ldots,v_{i_j}$, in which $$\bar{x}_{i_j}^l = x_{i_j}^l\cos\theta_{i} - y_{i_j}^l\sin\theta_{i} + \bar{x}_{i}^1,$$ and $$\bar{y}_{i_j}^l=x_{i_j}^l\sin\theta_{i} + y_{i_j}^l\cos\theta_{i} + \bar{y}_{i}^1,$$
\noindent being $(\bar{x}_i^1,\bar{y}_i^1)$ the variable reference point of polygon $P_i$ and $\theta_i$ the variable rotation angle.  $(x_{i_j}^l,y_{i_j}^l)$ are the coordinates of a vertex of a polygon $P_i$ in the original position.

Constraints \ref{ec1} and \ref{ec2} ensure that a polygon $P_i$ is integrally inside the strip.  In these constraints, the width $e$, is a fixed parameter; the length $z$ is a variable; $\bar{x}_{i_j}^l$ and  $\bar{y}_{i_j}^l$ depend on the reference point and the rotation angle of the polygon, $(\bar{x}_i^1,\bar{y}_i^1)$ and $\theta_i$, which are variable, respectively.
Constraints \ref{ec3} and \ref{ec4} ensure non-overlapping of the convex polygons $P_{i_j}$ and $P_{r_s}$. In these constraints, $c_{i_j,r_s}$ and $d_{i_j,r_s}$ depend on the reference point and the rotation angle of the separation line, $(\bar{x}_{i_j,r_s}, \bar{y}_{i_j,r_s})$ and $\alpha_{i_j,r_s}$, which are variable, respectively. 


\section{Computational experiments and results}\label{resultados}

All numerical experiments were performed on an Intel Core I7-4510U CPU @ 2.1GHz processor and 8 GB of memory.  We used a code for nonlinear programming to solve the problem, IPOPT \cite{ipopt} (an algorithm of interior points type), which is part of the COIN-OR \cite{coin}. 

IPOPT is the implementation of a barrier or interior points method for nonlinear optimization on a large scale that solves problems of the type 

\begin{equation*}
\begin{aligned}
&\text{Minimize} && &&f(v)\\ 
&\text{subject to} && g^{L}\leq&&g(v)\leq g^{U},\\
& && v^{L}\leq &&v\leq v^{U},
\\
\end{aligned}
\end{equation*}

\noindent being  $v \in \mathbb R^{n}$, and $f$, and $g$ continuously differentiable; the mathematical details of this algorithm can be found in \cite{ipopt,Nocedal}.  We used the C version of Ipopt-3.12.3 and compile the codes in the Ubuntu 12.04 operating system.

The CPU time is very large when we use the Hessian of the Lagrangian, therefore, we will always use an option given in IPOPT to approximate the Hessian with limited memory, which makes the runtime shorter, without affecting the quality of the solution.  In addition to the standard IPOPT parameters, we use the adaptive update strategy for barrier parameter, and we add equality constraints when handling fixed variables. The maximum execution time is set to one hour.

In the next subsection a brief explanation of the starting point used in the execution of IPOPT is presented.  In Section \ref{results}, the results obtained with IPOPT and comparisons of these with two methodologies recent in the literature (\cite{stoyan15,Liao16}), which also allow free rotations, are presented.


\subsection{Starting points}\label{starting point}

The solution given by IPOPT is a stationary point.  Taking into account that the developed nonlinear models are no convex, when solving them we can find many stationary points with different objective function values, and these stationary points depend on the starting point.  

To generate a starting point, we used a bottom-left algorithm, which is a single pass heuristic that, given a set of pieces and an order, places the pieces one by one on the strip, as far to left and to the bottom as possible.  We choose the shortest length obtained from the execution of 1000 iterations. In each one of these iterations the algorithm receives a list of randomly ordered polygons.  This list is represented by a sequence, and for decoding the sequences, we used the technique presented in \cite{mundim17}.  To avoid overlapping, the algorithm uses no-fit raster, concept also introduced in \cite{mundim17}. In no-fit raster, the strip should be represented using a discrete grid of points, the scale used for discretization in most instances is 1.0, except for Albano, Dagli, and Swim, which are 0.02, 0.5, and 0.00005, respectively. In all instances, the allowed rotation of the polygons are at four predefined angles $0^{\circ}, 90^{\circ}, 180^{\circ}$, and $270^{\circ}$. 

 
\subsection{Comparing results}\label{results}

To test our model, we will use the same benchmark problems used in the two approaches with which we will compare results, and that can be downloaded in \cite{esicup}.  The most important characteristics of these instances are presented in Table~\ref{instances}.  The names of the instances are presented in the first column.  In the second and third columns the number of convex, and no convex polygons are presented, respectively.  The total number of polygons, after decomposition of no convex polygons into convex polygons are presented in the fourth column. In the fifth the total number of vertices are presented.  The number of variables, and the number of constraints are presented in the sixth and seventh columns, respectively.






\begin{table}[!htb]
\caption{Data of the instances}
\begin{center}
\begin{tabular}{c|r|r|r|r|r|r}
\hline
Instance  &  \parbox[t]{1.5cm}{convex polygons} &  \parbox[t]{1.5cm}{no convex polygons} & total & vertices & variables & constraints \\\hline
albano   & 10   & 14 & 52 & 220   & 3907 & 11700\\
blaz     & 16  & 12 & 48 & 216   & 3385 & 10792 \\
dagli    & 21   & 9 & 51 & 228   & 3790 & 11988 \\
jakobs1  & 15  & 10& 35 & 146  & 1831 & 5468        \\
jakobs2  & 14  & 11& 42 & 158  & 2590 & 7367        \\
marques  & 10   & 14 & 50 & 214   &  3628 & 11002     \\
poly1a   & 10  & 5 & 22 & 81   &  712 & 1968          \\
poly2a   & 20  & 10 & 44 & 162   &  2875 & 7500      \\
poly3a   & 30  & 15 & 66 & 243   &  6490 & 16596      \\
poly4a   & 40  & 20 & 88 & 324   &  11557 &  29256   \\
poly5a   & 50  & 25 & 110 & 405   & 18076  & 45480   \\
poly10a   & 100  & 50 & 220 & 810   & 72451 & 180060 \\
poly20a   & 200  & 100 & 440 & 1620 & 290101 & 716520\\
shirts   & 60  & 39 & 169 & 739   &  42583  & 126236 \\
swim   &   6  & 42 & 291 & 1446   &  123787 & 415845 \\
trousers   & 48  & 16 & 104 & 468   &   16045 & 49476\\\hline
\end{tabular}
\end{center}
\label{instances}
\end{table}



For each instance, we obtained 10 starting points using the bottom-left algorithm, mentioned in the previous section.  We execute IPOPT to solve our model with each one of these 10 different starting points.  In Table \ref{ZIPOPT}, the minimum (second column), the average (third column), and the maximum (fourth column) strip length obtained are presented.  This is obtained when we solve the model in Section \ref{modelo} for the instances of Table \ref{instances} using the 10 different starting points. The average time used to construct the starting points (fifth column) and the average time (sixth column) that was spent to solve the instances are also presented in the Table \ref{ZIPOPT}, as well as those results of the recent literature that allow free rotations. For these methods, the strip length and time reported in \cite{stoyan15} are in seventh and eighth columns, and strip length and time reported in \cite{Liao16} are in ninth and tenth columns, respectively.



\begin{table}[!htb]
\caption{Comparison of our results to those in \cite{stoyan15} and \cite{Liao16}.}
\centering
\resizebox{13cm}{!}{
\begin{tabular}{c|r|r|r|r|r|r|r|r|r}
\hline
         & \multicolumn{5}{|c|}{Our approach} &\multicolumn{2}{|c|}{Best results in \cite{stoyan15}} &\multicolumn{2}{|c}{Best results in \cite{Liao16}} \\\hline
Instance & Min. Sol. & Avg. Sol.  & Max. Sol. & \parbox[t]{1.3cm} {SP Avg. time(s)} & \parbox[t]{1.3cm} {IPOPT Avg. time(s)} &   Length        & Time(s)&  Length    & Time(s)          \\\hline
albano     &  10355.80  & 10601.03   & 10849.99   & 158.68  &   178.56  &  &  & \bf{10032.24} & 124.39\\
blaz       &     27.82  &    29.02   &    30.44   &   1.37  &    49.52  & \bf{25.41} & 25.42 & 28.27 & 56.86\\
dagli      &     60.60  &    61.96   &    63.36   & 404.13  &   135.10  & \bf{56.90} & 139.00 & 59.24 & 132.58\\
jakobs1    &     \bf{12.99}  &    13.49   &    14.00   &   4.68  &    19.30  &  &  & 13.19 & 48.46\\
jakobs2    &     26.00  &    27.29   &    30.00   &  17.60  &    18.41  &  &  & \bf{24.25} & 53.67\\
marques    &     \bf{84.65}  &    86.68   &    91.00   &  66.74  &    58.11  &  &  & 84.93 & 118.12 \\
poly1a     &     14.01  &    14.53   &    15.01   &   6.23  &    19.07  &  &  & \bf{13.90} & 27.59\\
poly2a     &     \bf{26.16}  &    27.30   &    27.92   &  10.68  &    87.53  &  &  & 26.67 & 61.23\\
poly3a     &     \bf{39.01}  &    40.29   &    44.50   &  15.44  &   773.32  &  &  & 39.48 & 149.66\\
poly4a     &     \bf{50.99}  &    52.50   &    53.72   &  20.59  &  1621.54  &  &  & 51.13 & 210.74 \\
poly5a     &     \bf{63.66}  &    65.05   &    66.31   &  26.21  &  1773.85  &  &  & 65.64 & 287.32\\
poly10a    &    127.16  &   129.46   &   140.01   &  58.25  &  2350.17  & \bf{126.29} & 618.80 & & \\
poly20a    &    254.86  &   268.06   &   294.81   & 151.38  &  3484.92  & \bf{251.04} & 1209.17 & & \\
shirts     &     \bf{62.19}  &    64.43   &    65.62   &  11.55  &  1808.73  &  &  & 65.06 & 340.89\\
swim       &   6011.93  &  6311.16   &  6526.38   & 449.85  &  3600.00  & \bf{5661.95} & 431.97 & & \\
trousers   &    \bf{249.35}  &   258.03   &   261.94   & 118.00  &   603.36  &  &  & 251.94 & 265.48\\\hline 
\end{tabular}
}
\label{ZIPOPT}
\end{table}

For those instances tested in \cite{stoyan15}, it can be observed that the length in the best solution obtained in this work is slightly greater than those reported there.  However, when comparing with \cite{Liao16} we can see that in most instances the minimum length obtained using our model is smaller.  In the few remaining, the length is very close. 
On the other hand, note that our model finds good solutions for problems with a large number of polygons, although in these, the number of variables and constraints, and therefore the computational time, grows drastically.
These solutions may not be optimal local due to the maximum execution time.

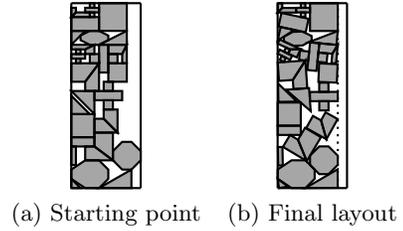
\begin{figure}[!h]
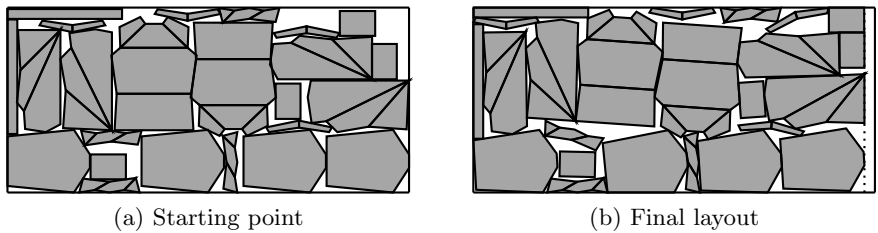

  \centering
  \caption{Instance: Albano}
  \begin{subfigure}{0.4\textwidth}
    \centering

  \caption{Starting point}
\end{subfigure}\hspace{0.05\textwidth}
\begin{subfigure}{0.4\textwidth}
  \centering
  %
\caption{Final layout}
\end{subfigure}
\label{fig:albano}
\end{figure}

\begin{figure}[!h]
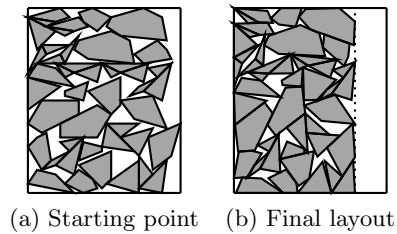

\centering
\caption{Instance: Poly3a}
\begin{subfigure}{0.25\textwidth}
  \centering
 %
  \caption{Starting point}
\end{subfigure}\hspace{0.05\textwidth}
\begin{subfigure}{0.25\textwidth}
  \centering
   %
  \caption{Final layout}
\end{subfigure}
\label{fig:poly3a}
\end{figure}

\begin{figure}[!h]
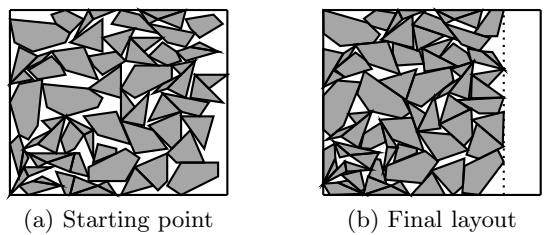

\centering
\caption{Instance: Poly4a}
\begin{subfigure}{0.3\textwidth}
  \centering
  %
  \caption{Final layout}
\end{subfigure}
\label{fig:poly4a}
\end{figure}

\begin{figure}[!h]
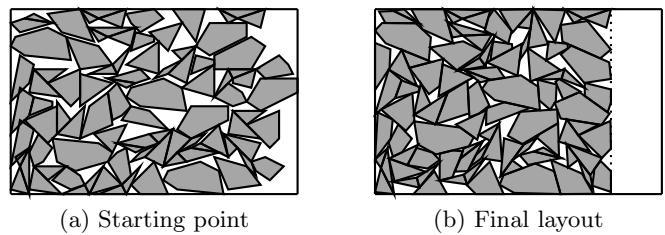

\centering
\caption{Instance: Poly5a}
\begin{subfigure}{0.4\textwidth}
  \centering
 %
  \caption{Final layout}
\end{subfigure}
\label{fig:poly5a}
\end{figure}

\begin{figure}[!h]
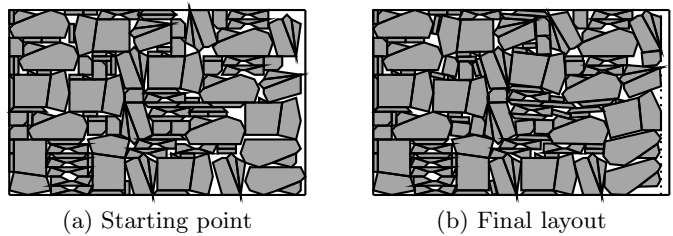

\centering
\caption{Instance: Swim}
\begin{subfigure}{0.2\textwidth}
  \centering
  %
  \caption{Final layout}
\end{subfigure}
\label{fig:swim}
\end{figure}

\begin{figure}[!h]
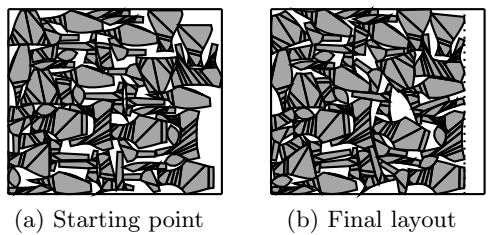

\caption{Instance: Trousers}
\begin{subfigure}[b]{\textwidth}
  \centering
  %
  \caption{Starting point}
\end{subfigure}
\begin{subfigure}[b]{\textwidth}
  \centering
   %
  \caption{Final layout}
\end{subfigure}
\label{fig:trousers}
\end{figure}

The layout of minimum length obtained and the starting point used to find it, for each instance, can be seen in Figures \ref{fig:albano} - \ref{fig:trousers}. 





\begin{figure}[!h]
\centering
\caption{Instance: Blaz}
\begin{subfigure}{0.35\textwidth}
  \centering
 %
  \caption{Starting point}
\end{subfigure}\hspace{0.05\textwidth}
\begin{subfigure}{0.35\textwidth}
  \centering
 %
  \caption{Final layout}
\end{subfigure}
\label{fig:blaz}
\end{figure}

\begin{figure}
\begin{minipage}[c]{7cm}
\centering
  \caption{Instance: Dagli}
  \begin{subfigure}{0.4\textwidth}
    \centering
 %
    \subcaption{Starting point}
  \end{subfigure}\hspace{0.05\textwidth}
  \begin{subfigure}{0.4\textwidth}
    \centering
%
  \subcaption{Final layout}
  \end{subfigure}
\end{minipage}
\begin{minipage}[c]{7cm}
  \centering
  \caption{Instance: Jakobs1}
  \begin{subfigure}{0.4\textwidth}
    \centering
%
    \subcaption{Starting point}
  \end{subfigure}
  \begin{subfigure}{0.35\textwidth}
    \centering
 %
    \subcaption{Final layout}
  \end{subfigure}
\end{minipage} 
\label{fig:jakobs2}
\end{figure}

\begin{figure}
\begin{minipage}[c]{7cm}
\centering
  \caption{Instance: Jakobs2}
  \begin{subfigure}{0.4\textwidth}
    \centering
%
    \subcaption{Starting point}
  \end{subfigure}
  \begin{subfigure}{0.4\textwidth}
    \centering
 %
  \subcaption{Final layout}
  \end{subfigure}
\end{minipage}
\begin{minipage}[c]{7cm}
  \centering
  \caption{Instance: Marques}
  \begin{subfigure}{0.4\textwidth}
    \centering
%
    \subcaption{Starting point}
  \end{subfigure}
  \begin{subfigure}{0.35\textwidth}
    \centering
 %
    \subcaption{Final layout}
  \end{subfigure}
\end{minipage} 
\label{fig:jakobs2}
\end{figure}

\begin{figure}
\begin{minipage}[c]{7cm}
\centering
  \caption{Instance: Poly1a}
  \begin{subfigure}{0.4\textwidth}
    \centering
%
  \subcaption{Final layout}
  \end{subfigure}
\end{minipage}
\begin{minipage}[c]{7cm}
  \centering
  \caption{Instance: Poly2a}
  \begin{subfigure}{0.4\textwidth}
    \centering
%
    \subcaption{Final layout}
  \end{subfigure}
\end{minipage} 
\label{fig:jakobs2}
\end{figure}

\begin{figure}[!h]
\caption{Instance: Poly10a}
\begin{subfigure}[b]{\textwidth}
  \centering
  %
  \caption{Starting point}
\end{subfigure}
\begin{subfigure}[b]{\textwidth}
  \centering
   %
  \caption{Final layout}
\end{subfigure}
\label{fig:poly10a}
\end{figure}

\begin{figure}[!h]
\caption{Instance: Poly20a}
\begin{subfigure}[b]{\textwidth}
  \centering
  %
  \caption{Starting point}
\end{subfigure}
\begin{subfigure}[b]{\textwidth}
  \centering
   %
  \caption{Final layout}
\end{subfigure}
\label{fig:poly20a}
\end{figure}

\begin{figure}[!h]
\centering
\caption{Instance: Shirts}
\begin{subfigure}{0.3\textwidth}
  \centering
  %
  \caption{Final layout}
\end{subfigure}
\label{fig:shirts}
\end{figure}

\section{Conclusions}\label{conclusion}

In this work, we developed a model for irregular strip packing problems, that allows free rotations and uses separation lines to avoid overlap.  As a relevant point of this work we highlight the use of the general equation of the line to model the separation lines, allowing us to use only three variables: the rectangular coordinates of the reference point and the angle of rotation. This establishes a marked dissimilarity with other models presented in the literature, for example the one presented in \cite{kallrath2009} which has a higher number of variables. In fact, in such work, it was reported that for one instance with a circle and a polygon of 4 vertices, 47 variables were used. In addition, it presented a shape of an experiment with two convex polygons of 5 and 6 vertices, but it does not report the number of variables of this instance.  We estimate that there are 87 variables. Instead, for an instance with two convex polygons our model would only need 10 variables, a significant simplification that leads to a better performance of the solution method ensuring likewise a good solution to the problem.  The variables in the model used in \cite{kallrath2009} are due to the vector equation of the line used to model the separation lines, which in turn implies the employment of many variables such as footing point vector, direction vector, normal vector, vectors connecting the separating line with the vertices, and distance of the vertices to the separating line, among others. Furthermore, the reference point and the angle of rotation are not the only relevant variables for each polygon, as it happens here, but also the vertices of the polygon are variables, which would increase the complexity of the model and therefore, of the solution method for the problem.

The solution of the problem modeled here, using local nonlinear programming solvers, depends on the starting point.  We use a bottom-left algorithm to construct these starting points. To test the effectiveness of our model, we compare our results with those obtained recently in the literature (\cite{stoyan15,Liao16}), which also use methodologies with free rotations.  The lengths reported in \cite{stoyan15} are smaller but very close to those found in this work.  On the other hand, the lengths reported in \cite{Liao16} are greater in most instances; in the others, they are very close. With the above said, the effectiveness of our model is verified, as well as the constructing starting points using bottom-left algorithm is good; however, we believe that these results could be improved by using another algorithm to construct the starting points, dividing the strip and the number of polygons into several parts, thus solving sub problems, among others.


\section*{Acknowledgements}

This research was partially supported by CNPq (grant 141072/2014-8 and 409043/2016-8) and FAPESP (grant 2013/07375-0 and 2016/01860-1), from Brazil.

\medskip

\begin{thebibliography}{99}


\bibitem[1]{b6} 
\newblock Albano A. \& Sapuppo A.(1980).
\newblock \emph{Optimal allocation of two-dimensional irregular shapes using heuristic search methods},
\newblock IEEE Transactions on Systems, Man and Cybernetics, \textbf{10}: 242--248.

\bibitem[2]{alvarez-valdes} 
\newblock Alvarez-Valdes R., Martinez A. \& Tamarit J.M.(2013).
\newblock \emph{A branch and bound algorithm for cutting and packing irregularly shaped pieces},
\newblock International Journal of Production Economics, \textbf{145}(2): 463--477.

\bibitem[3]{andreani2007} 
\newblock Andreani R., Birgin E.G., Martinez J.M. \& Schuverdt M.L.(2007).
\newblock \emph{On augmented lagrangian methods with general lower-level constraints},
\newblock SIAM Journal on Optimization, \textbf{18}: 1286--1309.

\bibitem[4]{andreani2008} 
\newblock Andreani R., Birgin E.G., Martinez J.M. \& Schuverdt M.L.(2008).
\newblock \emph{Augmented lagrangian methods under the Constant Positive Linear Dependence constraint qualification},
\newblock Mathematical Programming, \textbf{111}: 5--32.

\bibitem[5]{bennell2008} 
\newblock Bennell J.A. \& Oliveira J.F (2008).
\newblock \emph{The geometry of nesting problems: A tutorial},
\newblock European Journal of Operational Research, \textbf{184}: 397--415.

\bibitem[6]{bennell2010} 
\newblock Bennell J.A. Scheithauer G., Stoyan Y. \& Romanova T. (2010).
\newblock \emph{Tools of mathematical modelling of arbitrary object packing problems},
\newblock Ann Oper Res, \textbf{179}: 343--368.

\bibitem[7]{Burke06} 
\newblock Burke E., Hellier R., Kendall G. \& Whitwell G. (2006).
\newblock \emph{A New Bottom-Left-Fill Heuristic Algorithm for the Two-Dimensional Irregular Packing Problem},
\newblock Operations Research, \textbf{54}(3): 587--601.


\bibitem[8]{chernov2010} 
\newblock Chernov N., Stoyan Y. \& Romanova T. (2010).
\newblock \emph{Mathematical model and efficient algorithms for object packing problem},
\newblock Computational Geometry: Theory and Applications, \textbf{43}: 535--553.

\bibitem[9]{b8} 
\newblock Egeblad J., Nielsen, B.K. \& Odgaard A. (2007).
\newblock \emph{Fast neighborhood search for two and three-dimensional nesting problems},
\newblock European Journal of Operational Research, \textbf{183}: 1294--1266.


\bibitem[10]{Elkeran2013} 
\newblock Elkeran A. (2013).
\newblock \emph{A new approach for sheet nesting problem using guided cuckoo search and pairwise clustering},
\newblock European Journal of Operational Research, \textbf{231}: 757--769.


\bibitem[11]{esicup} 
\newblock EURO special interest group on cutting and packing. (2015). 
\newblock Available at
                \url{http://paginas.fe.up.pt/~esicup/datasets}

\bibitem[12]{fischetti} 
\newblock Fischetti M. \& Luzzi I. (2009).
\newblock \emph{Mixed-integer programming models for nesting problems},
\newblock Journal of Heuristics, \textbf{15}(3): 201--226.


\bibitem[13]{b3} 
\newblock Gomes A.M. \& Oliveira J.F. (2002).
\newblock \emph{A 2-exchange heuristic for nesting problems},
\newblock European Journal of Operational Research, \textbf{141}: 359--370.

\bibitem[14]{b4} 
\newblock Gomes A.M. \& Oliveira J.F. (2006).
\newblock \emph{Solving irregular strip packing problems by hybridising simulated annealing and linear programming.},
\newblock European Journal of Operational Research, \textbf{171}: 811--829.


\bibitem[15]{jones2013} 
\newblock Jones D.R. (2013).
\newblock \emph{A fully general, exact algorithm for nesting irregular shapes
},
\newblock Journal of Global Optimization, \textbf{59}: 367--404.

\bibitem[16]{Liao16} 
\newblock Liao X., Ma J., Ou C., Long F. \& Liu X. (2016).
\newblock \emph{Visual nesting system for irregular cutting-stock problem based on rubber band packing algorithm.},
\newblock Advances in Mechanical Engineering, \textbf{8}(6): 1--15.

\bibitem[17]{kallrath2009} 
\newblock Kallrath J. (2009).
\newblock \emph{Cutting circles and polygons from area-minimizing rectangles},
\newblock Journal of Global Optimization, \textbf{43}: 299--328.

\bibitem[18]{b7} 
\newblock Marques V.M., Bispo C.F. \& Sentieiro J.J. (1991).
\newblock \emph{A system for the compaction of two-dimensional irregular shapes based on simulated annealing},
\newblock IEEE Transactions on Industrial Electronics, Control and Instrumentation, \textbf{3}: 1911--1916.


\bibitem[19]{MIQO} 
\newblock Misener R. \& Floudas C.A. (2013).
\newblock \emph{GloMIQO: global mixed-integer quadratic optimizer},
\newblock Journal Global Optimization, \textbf{57}: 3--50.


\bibitem[20]{mundim17}
\newblock Mundim L. R., Andretta M. \& Queiroz T.A.(2017)
\newblock \emph{A biased random key genetic algorithm for open dimension nesting problems using no-fit raster},
\newblock Expert Systems with Applications, \textbf{81}: 358--371.

\bibitem[21]{b9} 
\newblock Nielsen B.K. (2007).
\newblock \emph{An efficient solution method for relaxed variants of the nesting problem},
\newblock Proceedings of the thirteenth Australasian symposium on Theory of computing, \textbf{65}: 123--130.


\bibitem[22]{Nocedal} 
\newblock Nocedal J., Wächter A \& Waltz R.A. (2009).
\newblock \emph{Adaptive barrier strategies for nonlinear interior
methods},
\newblock SIAM Journal on Optimization \textbf{19}: 1674--1693.

\bibitem[23]{b1} 
\newblock Oliveira J.F. \& Ferreira J.S. (2008).
\newblock \emph{Algorithms for nesting problems, applied simulated annealing
},
\newblock In: Vidal, R.V.V. (Ed.), Lecture notes in econ. and Maths Systems. Springer Verlag, \textbf{396}: 255--274.


\bibitem[24]{rochaetal} 
\newblock Rocha P., Rodrigues R. Gomes A.M., Toledo F.M.B \& Andretta M. (2015).
\newblock \emph{Two-Phase Approach to the Nesting problem with continuous rotations},
\newblock IFAC-PapersOnline, \textbf{48}(3): 501--506.

\bibitem[25]{sahinidis:baron:14.3.1} 
\newblock Sahinidis N.V. (2014).
\newblock \emph{BARON 14.3.1: Global Optimization of Mixed-Integer Nonlinear Programs},
\newblock User's Manual, Available at
                \url{http://www.minlp.com/downloads/docs/baron\%20manual.pdf}

\bibitem[26]{SB} 
\newblock Segenreich S.A. \& Braga L.M. (1986).
\newblock \emph{Optimal nesting of general plane figures: a Monte Carlo heuristical approach},
\newblock Computers and Graphics, \textbf{10}: 229--237.

\bibitem[27]{stoyan2001} 
\newblock Stoyan Y.G., Terno J., Scheithauer, G., Gil N., \& Romanova T. (2001).
\newblock \emph{Phi-functions for primary 2d-objects},
\newblock Studia Informatica Universalis, \textbf{2}(1): 1--32.

\bibitem[28]{stoyan2004} 
\newblock Stoyan Y.G., Scheithauer, G., Gil N., \& Romanova T. (2004).
\newblock \emph{Phi-functions for complex 2d-objects},
\newblock 4OR: Quartely Journal of the Belgian, French and Italian Operations Research Societies, \textbf{2}: 69--84.

\bibitem[29]{Stoyan2008} 
\newblock Stoyan Y.G. \& Chugay A.M. (2008).
\newblock \emph{Packing cylinders and rectangular parallelepipeds with distances between them},
\newblock European Journal Operation Research, \textbf{197}: 446--455.

\bibitem[30]{stoyan2012} 
\newblock Stoyan Y.G., Zlotnik M.V. \& Chugay A.M. (2012).
\newblock \emph{Solving an optimization packing problem of circles and non-convex polygons with rotations into a multiply connected region},
\newblock  Journal of the Operational Research Society, \textbf{63}: 379--391.

\bibitem[31]{stoyan15} 
\newblock Stoyan Y.G., Pankratov A. \& Romanova T. (2016).
\newblock \emph{Cutting and packing problems for irregular objects with continuous rotations: mathematical modelling and non-linear optimization.},
\newblock Journal of the Operational Research Society, \textbf{67}(5): 786--800.

\bibitem[32]{ts05} 
\newblock Tawarmalani M. \& Sahinidis N.V. (2005).
\newblock \emph{A polyhedral brach-and-cut approach to global optimization},
\newblock Mathematical Programming, \textbf{103}(2): 225--249.


\bibitem[33]{toledo} 
\newblock Toledo F.M.B, Carravilla M.A, Ribeiro C., Oliveira J.F. \& Gomes A.M. (2013).
\newblock \emph{The dotted-board model: A new mip model for nesting irregular shapes},
\newblock International Journal of Production Economics, \textbf{145}(2): 478--487.


\bibitem[34]{ipopt} 
\newblock W\"achter A. \& Biegler L.T. (2006).
\newblock \emph{On the implementation of a primal-dual interior point filter line search algorithm for large-scale nonlinear programming},
\newblock Mathematical Programming, \textbf{106}(1): 25--57.


\bibitem[35]{coin} 
\newblock W\"achter A. \& Biegler L.T. (2015).
\newblock \emph{COIN OR project}, Available at
                \url{http://projects.coin-or.org/Ipopt}



\bibitem[36]{Z1960} 
\newblock Zoutendijk G. (1960).
\newblock \emph{Methods of feasible directions, a study in linear and non-linear programming},
\newblock Elsevier.

\bibitem[37]{Z1970} 
\newblock Zoutendijk G. (1970).
\newblock \emph{Nonlinear programming, computational methods},
\newblock Integer and Nonlinear Programming, \textbf{143}(1): 37--86.



\end{thebibliography}
\end{document}